\documentclass[a4paper, 12pt]{article}

\usepackage{amssymb}
\usepackage{amsbsy}
\usepackage{graphicx}

\begin{document}

\title{On contractive families and a fixed-point question of Stein}
\author{Tim D. Austin}
\date{}

\maketitle


\newenvironment{nmath}{\begin{center}\begin{math}}{\end{math}\end{center}}

\newtheorem{thm}{Theorem}[section]
\newtheorem{lem}[thm]{Lemma}
\newtheorem{prop}[thm]{Proposition}
\newtheorem{cor}[thm]{Corollary}
\newtheorem{conj}[thm]{Conjecture}
\newtheorem{dfn}[thm]{Definition}


\newcommand{\Bor}{\mathcal{B}}
\renewcommand{\Pr}{\mathrm{Pr}}
\newcommand{\s}{\sigma}
\renewcommand{\P}{\mathcal{P}}
\renewcommand{\O}{\Omega}
\renewcommand{\S}{\Sigma}
\newcommand{\T}{\mathrm{T}}
\newcommand{\co}{\mathrm{co}}
\newcommand{\e}{\mathrm{e}}
\renewcommand{\i}{\mathrm{i}}
\renewcommand{\l}{\lambda}

\newcommand{\G}{\Gamma}
\newcommand{\g}{\gamma}
\renewcommand{\L}{\Lambda}
\newcommand{\hcf}{\mathrm{hcf}}
\newcommand{\F}{\mathcal{F}}

\newcommand{\bb}[1]{\mathbb{#1}}
\renewcommand{\rm}[1]{\mathrm{#1}}
\renewcommand{\cal}[1]{\mathcal{#1}}

\newcommand{\qed}{\nolinebreak$\Box$}

\begin{quote}
\noindent \textbf{Abstract} \hspace{3pt} The following conjecture
generalizing the Contraction Mapping Theorem was made by Stein in
\cite{Stein}:

\emph{Let $(X,\rho)$ be a complete metric space and let $\F =
\{T_1,\ldots,T_n\}$ be a finite family of self-maps of $X$.
Suppose there is a constant $\g \in (0,1)$ such that for any $x,y
\in X$ there exists $T \in \F$ with $\rho(T(x),T(y)) \leq \g
\rho(x,y)$.  Then some composition of members of $\F$ has a fixed
point.}

In this paper we disprove this conjecture.  We also show that it
does hold for a (continuous) commuting $\F$ in the case $n=2$. We
conjecture that it holds for commuting $\F$ for any $n$.
\end{quote}

\section{Introduction}

Given a metric space $(X,\rho)$, a \textbf{self-map} of $X$ will
be a continuous function $T$ from $X$ to itself.  Such a $T$ is
said to be a \textbf{contraction} if there is some $\g < 1$ such
that $\rho(T(x),T(y)) \leq \g\rho(x,y)$ for all $x,y \in X$.

The original contraction mapping principle is due to Banach:

\begin{prop}
Let $(X,\rho)$ be a complete metric space and $T:X \to X$ a
self-map for which there exists a constant $\g \in (0,1)$ such
that for any $x,y \in X$ we have $\rho(T(x),T(y)) \leq
\g\rho(x,y)$.  Then $T$ has a unique fixed point in $X$.
\end{prop}

Many possible generalizations of this result have been explored
since the appearance of Banach's original version (\cite{Ban}).
There is now a major body of fixed point theory depending on
algebraic topology (and thus usually applicable only in spaces
with much ``nicer'' structure than an arbitrary complete metric
space), and also a range of deep results in spaces with other
underlying properties (such as Brouwer's theorem for certain
subsets of $\bb{R}^n$ and its extension to locally convex
Hausdorff topological vector spaces due to Schauder; see, for
example, \cite{Hat} for Brouwer's theorem and \cite{Con} for
Schauder's theorem).

One elegant generalization is the Generalized Banach Contraction
Theorem, conjectured by Jachymski, Schroder and Stein in
\cite{JachSchroStein} and proved by Merryfield, Rothschild and
Stein in \cite{MerryRothStein}:

\begin{prop}
Suppose $(X,\rho)$ is a complete metric space and $T$ is a
self-map of $X$. Suppose further that there exist $m \in \bb{N}$,
$\g \in (0,1)$ such that for any $x,y \in X$ we have
\[\min\{\rho(T^k(x),T^k(y)):\ 1 \leq k \leq m\} \leq \g\rho(x,y).\]
Then $T$ has a (necessarily unique) fixed point in $X$.
\end{prop}

\noindent \textbf{Remark} \hspace{3pt} This result is now known to
hold for discontinuous $T$ for which such $\g$ and $m$ exist; this
has been proved by Merryfield and Stein in \cite{MerryStein}.

\vspace{10pt}

The arguments used in these papers rest crucially on a subtle
application of the infinite case of Ramsey's Theorem.  At about
the same time as the last of these papers appeared, a different
approach was discovered by Arvanitakis, and has now been published
in \cite{Arvan}; this is also highly combinatorial.

The following further conjecture was made in \cite{Stein}:

\begin{conj}
Let $(X,\rho)$ be a complete metric space and let $\F =
\{T_1,\ldots,T_n\}$ be a finite family of self-maps of $X$.
Suppose there is a constant $\g \in (0,1)$ such that for any $x,y
\in X$ we may choose $U \in \F$ with $\rho(U(x),U(y)) \leq \g
\rho(x,y)$.  Then some composition of members of $\F$ has a fixed
point.
\end{conj}

We will refer to a family $\F$ of self-maps satisfying the
hypotheses of the above conjecture as a \textbf{$\g$-contractive
family}; $\F$ is \textbf{contractive} if it is $\g$-contractive
for some $\g$.

We show in Section \ref{sec:counterex} below that the conjectured
result is false in general. We construct an example of two
self-maps $S$ and $T$ of a complete metric space $X$ which satisfy
the conditions of the conjecture with $\g = \frac{3}{4}$, but such
that no member of the semigroup generated by $\{S,T\}$ has a fixed
point. In Section \ref{sec:commcaseprog} we show that the result
is true for two self-maps if they commute, and conjecture that it
holds for any finite commuting family of self-maps.

\section{Counterexample for general result}\label{sec:counterex}

Here we construct an example of a complete non-compact metric
space $(X,\rho)$ and a pair of Lipschitz-1 self-maps $S$ and $T$
of $X$ such that $\{S,T\}$ is contractive but no word in $S$ or
$T$ has a fixed point.  We do not know whether an example of such
behaviour exists with $X$ compact.

The construction to be described here is motivated by the problem
itself. We will construct a metric on a certain non-complete space
and family of self-maps all at once, then complete the space and
analyse it carefully to show that there are no fixed points. This
construction avoids the structural features which force many of
the more familiar types of self-map on simpler spaces to have
fixed points.

\subsection{The basic construction}

We begin with an example of a set with a pair of non-commuting
self-maps.  Write $\G_0$ for the free monoid on two generators,
say $a$ and $b$; as a set, $\G_0$ contains all words in $a$ and
$b$ of finite length, and the empty word, which we write as $0$
for convenience.

We recall some standard definitions.  Given a word $w \in \G_0$,
we will write $l(w)$ for its length (later we will define a
different notion of length, the $d$-length). Given $u,v \in \G_0$
we say that $v$ \textbf{extends} $u$ and write $u\ |\ v$ if $v$ is
an extension of $u$, that is if $v = uw$ for some $w \in \G_0$. We
define the \textbf{meet} of $u$ and $v$ to be the longest initial
segment of $u$ that is also an initial segment of $v$, and write
it as $u \wedge v$. Given $w \in \G_0$ and $n \geq 0$, we let
$w|_n$ be the \textbf{restriction of $w$ to $n$}: if $w =
c_1c_2\cdots c_m$ with each $c_1$ equal to either $a$ or $b$ then
\[w|_n = \left\{\begin{array}{ll}
0\ \ \ &\hbox{if $n=0$}\\ c_1c_2\cdots c_n\ \ \ &\hbox{if $1 \leq
n \leq m$}\\ w\ \ \ &\hbox{if $n > m$}\end{array}\right.\]

We may define a graph $G = (V,E)$ by taking $V = \G_0$ for the set
of vertices and saying that for $u,v \in \G_0$ we have $\{u,v\}
\in E$ if and only if one of $u,v$ is an extension of the other by
a single letter; that is, if and only if one of the following
holds:
\[u = va,\ \ \ u = vb,\ \ \ v=ua,\ \ \ v=ub.\]
This graph is an infinite rooted tree with root the empty word in
$\G_0$ and all degrees equal to two, as sketched in Fig 1.

\begin{figure}
\begin{center}
\includegraphics[width = 0.5\textwidth]{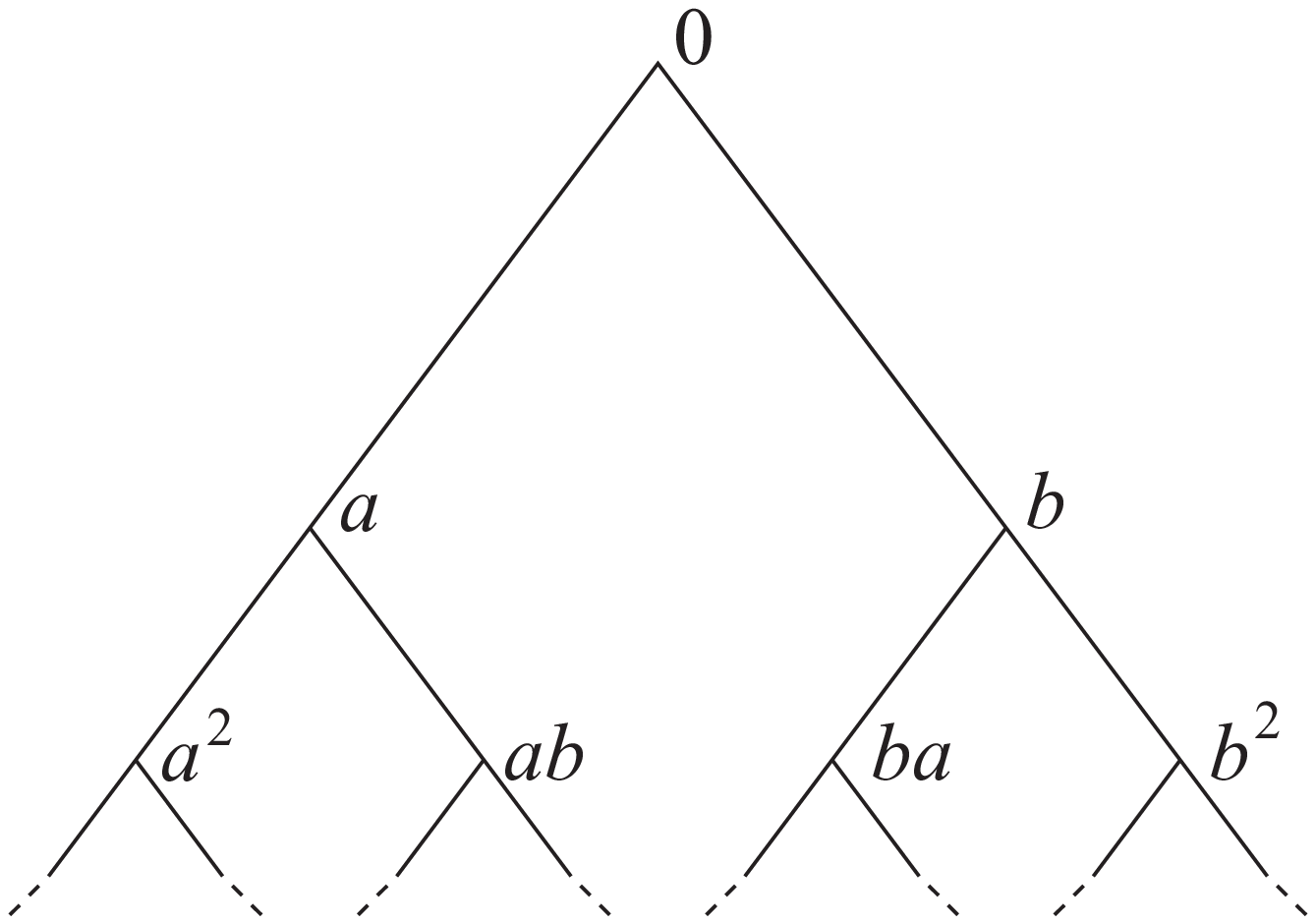}
\caption{$\G_0$}
\end{center}
\end{figure}

Now define two self-maps $S_0,T_0$ of $\G_0$ as pre-multiplication
by $a$ and $b$ respectively:
\[S_0(w) = aw,\ \ \ T_0(w) = bw\]
for $w \in \G_0$. It is clear that $S_0$ and $T_0$ do not commute;
indeed, we have a canonical identification of the semigroup
generated by $S_0$ and $T_0$ with $\G_0$ via the map $S_0 \mapsto
a$, $T_0 \mapsto b$.

We will construct our example as the completion of $\G_0$ under a
suitably-chosen metric, with $S_0$ and $T_0$ extended to maps $S$
and $T$ by continuity.  However, first we require some additional
notation and some routine groundwork.

Suppose $d$ is a function $E \to (0,\infty)$ (later we will pick a
particular $d$).  We may interpret $d$ as assigning a length to
each edge of $G$; for sake of brevity we write $d(u,v)$ in place
of $d(\{u,v\})$ for $u,v \in \G_0$ with $\{u,v\} \in E$. Since $G$
is connected $d$ now gives rise to a metric $\rho$ on $\G_0$
defined by
\[\rho(u,v) = \inf\left\{\sum_{i=0}^{n-1}d(u_i,u_{i+1}):\
(u_0,u_1,\ldots,u_n)\ \hbox{is a walk in $G$ with}\ u_0 = u,u_n =
v\right\}.\] This makes $(\G_0,\rho)$ into a metric space in which
$\rho$ is an example of a \textbf{tree metric}.  These are of some
interest in discrete mathematics and computer science; see, for
example, \cite{Dah} and the references listed there.

The \textbf{$d$-length} $l_d(w)$ of a word $w \in \G_0$ is given
by
\[l_d(w) = \rho(0,w) = \sum_{i = 0}^{l(w)-1}d(w|_i,w|_{i+1}),\]
the sum of the $d$-lengths of the edges traversed along the unique
path from $0$ to $w$.  It is clear from the above definitions that
for any $u,v \in \G_0$ we have
\[\rho(u,v) = \rho(u,u \wedge v) + \rho(u \wedge v,v) = l_d(u) + l_d(v)
- 2l_d(u \wedge v).\]

We will now develop some theory of the structure of the completion
of $\G_0$ under $\rho$.  We will write this as $\G$ and will use
the notation $\rho$ also for the extension of the metric to $\G$.
As is standard, $\G$ can be identified as the set of all
equivalence classes of $\rho$-Cauchy sequences in $\G_0$ under the
equivalence relation $\sim$, where $(u_n)_{n \geq 0} \sim (v_n)_{n
\geq 0}$ for $\rho$-Cauchy sequences $(u_n)_{n \geq 0}, (v_n)_{n
\geq 0}$ if and only if $\rho(u_n,v_n) \to 0$ as $n \to \infty$.

Write $\L$ for the set of all infinite words in $a$ and $b$; that
is, the set of all $\bb{N}\,$-indexed sequences of $a$s and $b$s.
In analogy with the above finitary definitions we say that $z \in
\L$ \textbf{extends} $u \in \G_0$ if $u$ is an initial segment of
$z$, and given such a sequence $z = (c_i)_{i \geq 1}$ and $n \geq
0$ we define the \textbf{restriction $z|_n$ of $z$ to $n$} by
\[z|_n = \left\{\begin{array}{ll}0\ \ \ &\hbox{if $n=0$}\\ c_1c_2\cdots
c_n\ \ \ &\hbox{if $n \geq 1$}\end{array}\right..\] For any $x,y
\in \G_0 \cup \L$ we may still define the meet $x \wedge y$ as the
largest initial segment of $x$ that is also an initial segment of
$y$; $x \wedge y$ is finite unless $x=y$ is infinite.  Also, we
define the \textbf{$d$-length} $l_d(z) \in [0,\infty]$ by
\[l_d(z) = \sum_{i\geq 0}d(z|_i,z|_{i+1}).\]

\begin{lem}\label{lem:convcond}
Suppose $z = (c_i)_{i \geq 1} \in \L$, and set $u_n = z|_n$ for
each $n \geq 0$.  Then $(u_n)_{n \geq 0}$ is a $\rho$-Cauchy
sequence if and only if
\[\sum_{n \geq 0}d(u_n,u_{n+1}) < \infty.\]
\end{lem}

\noindent \textbf{Proof} \hspace{3pt} Clearly we have $u_n\ |\
u_m$ whenever $m > n$, and hence
\[\rho(u_n,u_m) = \sum_{i=n}^{m-1}d(u_i,u_{i+1}).\]
Thus the condition that $(u_n)_{n \geq 0}$ be $\rho$-Cauchy
translates into an equivalent condition about the tails of the
infinite series $\sum_{n \geq 0}d(u_n,u_{n+1})$, and it is now
clear that this is equivalent to the convergence of the series.
\qed

\vspace{10pt}

The following result is perhaps less obvious.

\begin{lem} \label{lem:Cauchychar}
Let $(v_n)_{n \geq 0}$ be a $\rho$-Cauchy sequence in $\G_0$. We
have the following dichotomy: either $(v_n)_{n \geq 0}$ is
eventually equal to a constant, say $v$, in which case $(v_n)_{n
\geq 0} \sim v$, or there is a unique $z = (c_i)_{i \geq 1} \in
\L$ such that if the sequence $(u_n)_{n \geq 0}$ is defined from
$z$ as in the previous lemma then $(v_n)_{n \geq 0} \sim (u_n)_{n
\geq 0}$.
\end{lem}

\noindent \textbf{Proof} \hspace{3pt} Suppose that $(v_n)_{n \geq
0}$ is not eventually constant.  Then we must have $l(v_n) \to
\infty$ as $n \to \infty$, for otherwise there would be some
finite word $v$ occurring in $(v_n)_{n \geq 0}$ infinitely often
(since there are only finitely many words of a given finite
length), in which case the sequence would converge to $v$ and so
must be eventually equal to $v$. Hence, by passing to a
subsequence if necessary, we may assume that the lengths $l(v_n)$
are strictly increasing with $n$.

Now it is clear that the initial segments $v_n|_i$ must be
eventually constant as $n$ increases.  We may therefore define $z
= (c_i)_{i \geq 1}$ by requiring that $z|_i$ be the word equal to
$v_n|_i$ for all sufficiently large $n$.  It follows that
\[l_d(z) = \sum_{i\geq 0}d(z|_i,z|_{i+1}) = \lim_{n \to \infty}
l_d(z|_n) < \infty,\] for if $l_d(z|_n)$ tended to $\infty$ with
$n$ then we would have $l_d(v_k) \geq l_d(v_k|_n) = l_d(z|_n) \to
\infty$ for all sufficiently large $k \geq n$, contradicting that
$(v_n)_{n \geq 0}$ is a Cauchy sequence.  Hence given $\varepsilon
> 0$ we may choose $N \geq 0$ such that $\rho(v_n,v_k) \leq
\varepsilon$ whenever $n,k \geq N$, and also
\[\sum_{i \geq N}d(z|_i,z|_{i+1}) = l_d(z) - l_d(z|_N) \leq \varepsilon.\]
Since $v_n \wedge z = v_n \wedge v_k$ for all sufficiently large
$k \geq n$ and we have $\rho(v_n,v_n \wedge v_k) \leq
\rho(v_n,v_k)$, we deduce that $\rho(v_n,z \wedge v_n) \to 0$ as
$n \to \infty$, from which it readily follows that $\rho(v_n,z|_n)
\to 0$, so $(v_n)_{n \geq 0} \sim (u_n)_{n \geq 0}$.

Finally we observe that if $z^\prime = (c_i^\prime)_{i \geq 1}$ is
a different infinite word, $z^\prime \not= z$, then $z \wedge
z^\prime$ is finite and so writing $u_n^\prime = z^\prime|_{n}$
for $n \geq 0$ we see from the above construction of $z$ that
$\rho(v_n,u_n^\prime) \to \rho(z,z \wedge z^\prime) > 0$ as $n \to
\infty$.  Hence the required $z$ is unique. \qed

\vspace{10pt}

\noindent \textbf{Remark} \noindent We see from the above that
there is a bijective correspondence between the set of
non-eventually-constant $\rho$-Cauchy sequences in $\G_0$ and the
set of infinite words $z$ in $\L$ with $\sum_{n \geq
0}d(z|_n,z|_{n+1}) < \infty$.  Writing $\L_d$ for the set of such
$z \in \L$ we obtain a canonical identification of $\G$ with $\G_0
\cup \L_d$, suitably topologized.

\vspace{10pt}

Concerning the continuous extension of $S_0,T_0$ to $\G$ we have
the following.

\begin{lem}\label{lem:extofST}
Suppose $S_0,T_0$ may be extended by continuity to self-maps $S,T$
of $\G = \G_0 \cup \L_d$. If the Cauchy sequence $(u_n)_{n \geq
0}$ corresponds to $z \in \L_d$ then $(S_0(u_n))_{n \geq 0}$ is
also a Cauchy sequence  and corresponds to $az$ (the infinite word
$z$ pre-multiplied by $a$). The analogous result holds for $T_0$.
\end{lem}

\noindent \textbf{Proof} \hspace{3pt} We give the proof for $S_0$;
the argument for $T_0$ is similar.  The result is immediate: the
sequence $(S_0(u_n))_{n \geq 0}$ must be Cauchy since we are
assuming that $S_0$ has a continuous extension to $\G$, and if
$u_n = z|_n$ then by definition $S_0(u_n) = az|_n = (az)|_n$, so
$S_0(u_n)$ corresponds to $S(z) = az$. \qed

\vspace{10pt}

\noindent \textbf{Remark} \hspace{3pt} It follows that we may
write $S(z) = az$.

\subsection{The role of periodicity}

We will obtain our example by a careful choice of $d$ in the next
section; this will be the heart of the construction. First we will
need some results about periodicity within finite and infinite
words. These considerations are motivated by the following
observations.

We can see that no finite word in $\G_0 \subseteq \G$ is fixed by
any word in $S_0$ and $T_0$. If our construction is to provide a
suitable example we will need to be able to extend $S_0$ and $T_0$
to $S$ and $T$ on $\G$ so that they still do not have any fixed
points. Assuming that these extensions exist, we observe from
Lemma \ref{lem:extofST} we that the Cauchy sequence corresponding
to the infinite word $z \in \L_d$ is fixed by the word $U = R_k
\circ R_{k-1} \circ \cdots \circ R_1$ in $S,T$ if and only if
\[U(z) = R_k \circ R_{k-1} \circ \cdots
\circ R_1(z) = c_kc_{k-1}\cdots c_1z = z,\] where each $R_i$
corresponds to pre-multiplication by $c_i \in \{a,b\}$.  This
latter condition holds if and only if $z$ is periodic, so in order
to prevent words in $S,T$ from having any fixed points in $\L_d$
we wish to construct a $d$ so that with the metric $\rho$
periodicity is expensive. Before we do this we will need some
basic results about periodicity.

If $u,v \in \G_0$ we say that $u$ is a \textbf{shift} of $v$ if $u
= wv$ for some $w \in \G_0$, so $u$ is given by $v$ pre-multiplied
by some other word.  If $z \in \G_0$ or $z \in \L$ we say that $z$
\textbf{absorbs} $v$ if some shift of $v$ is an initial segment of
$z$; that is if some copy of $v$ appears as a segment of $z$. If
$z$ is periodic, we have the natural notions of its \textbf{repeat
block} (the smallest finite block of letters which can be repeated
to obtain $z$) and \textbf{period} $p(z)$ (the length of its
repeat block). Given $v \in \G_0$ we define its \textbf{minimal
potential period} $p_0(v)$ to be the shortest period of any $z$
which absorbs $v$ and its \textbf{minimal potential repeat blocks}
to be the repeat blocks of such $z$ (we show below that this is
well-defined). If $v = 0$ we set $p_0(v) = 0$. Clearly $p_0(v)
\leq l(v)$, since the infinite word obtained by repeating $v$
surely absorbs $v$.

\begin{lem}
The minimal potential repeat blocks of $v \in \G_0$ differ only by
cyclic permutations.
\end{lem}

\noindent \textbf{Proof} \hspace{3pt}  Clearly there are infinite
periodic words which absorb $v$ and whose repeat blocks appear in
full as a segment of $v$, so we may choose $z$ to be one such of
minimal period $p(z) = p_0(v)$ and repeat block $r$. Suppose $y$
is another such infinite periodic word.  Then $v$ appears as a
segment in both $z$ and $y$, and so the same is true of $r$. Since
$p(z) = p(y) = l(r)$, we infer that both $z$ and $y$ are of the
form $wrrr\cdots$ for $w$ some finite word in $\G_0$; hence any
repeat block taken from $y$ must be a cyclic permutation of $r$,
as required. \qed

\begin{lem}\label{lem:translates}
If $v \in \G_0\setminus \{0\}$ and $z \in \L$ is an infinite
periodic word of minimal period $p(z) = p_0(v)$ absorbing $v$ with
repeat block $r$, then any two segments of $z$ equal to $v$ differ
in position by a multiple of $l(r)$.
\end{lem}

By this we mean that if we have two segments of length $l(v)$ in
$z = (c_i)_{i \geq 1}$ equal to copies of $v$, say
$(c_i)_{i=n}^{n+l(v)-1}$ and $(c_j)_{j = k}^{k+l(v)-1}$ with $k >
n$, then $k-n$ is a multiple of $l(r)$.

\vspace{10pt}

\noindent \textbf{Proof} \hspace{3pt} The proof is by
contradiction.  Clearly any translate of a copy of $v$ by a
multiple of $l(r)$ will be another copy of $v$, since $z$ has
period $p(z) = l(r)$.  So we may suppose for sake of contradiction
that $n < k < n + l(r)$.  We know that some copy of $r$ appears as
a segment in $v$, so these two copies of $v$ yield two distinct
copies of $r$ within $z$ which differ in position by $k-n < l(r)$.
But since $z$ is an infinite string of repeated copies of $r$ this
means that the period of $z$ is at most $k-n < l(r)$,
contradicting the definition of $r$.  This completes the proof.
\qed

\begin{lem}\label{lem:perprop}
Suppose $v \in \G_0$ has a minimal potential repeat block $r$ and
$w \in \G_0$. Then:
\begin{enumerate}
\item $p_0(wv) \geq p_0(v)$;
\item precisely one of $av,bv$ has $r$ as a minimal potential repeat
block and so has minimal potential period equal to $p_0(v)$, and
the other has minimal potential period at least $p_0(v) + 1$.
\end{enumerate}
\end{lem}

\noindent \textbf{Proof} \hspace{3pt} \begin{enumerate} \item The
point is just that any infinite periodic word absorbing $wv$ must
absorb $v$;
\item Let $r$ be a minimal potential repeat block for $v$.  Then
$r$ appears as a segment of $v$, and $v$ appears as a segments of
any sufficiently large concatenation (``string'') of $r$'s. Taking
$z$ to be an infinite such string, we can ensure that there is a
copy of $v$ appearing as a segment within it which does not
contain the first letter of $z$.  Let $c$ be the letter ($a$ or
$b$) appearing in $z$ immediately before this copy of $v$; then
clearly $z$ absorbs $cv$, so $p_0(cv) = p_0(v)$. It remains to
show that if $d$ is the other letter (that is, the unique member
of $\{a,b\}\setminus \{c\}$) then $p_0(dv) > p_0(v)$.  Suppose
not.  Then we would have some infinite periodic $z_0$ of period
$l(r) = p_0(v)$ containing a segment equal to $v$ and immediately
preceded by the letter $d$.  However, since $r$ appears as a
segment of $v$, it would appear as a segment of $z_0$, and now
since $z_0$ is periodic of period $r$ we must have that $z_0$ is
of the form $wrrr\cdots = wz$ for some finite word $w$.  We know
that there is a copy of $v$ in $z$ immediately preceded by the
letter $c$.  Therefore $z_0$ contains some copies of $v$
immediately preceded by $c$ and some by $d$.  However, two such
copies cannot differ in position by a multiple of $l(r)$ by the
periodicity of $z_0$.  This contradicts Lemma
\ref{lem:translates}, and so completes the proof. \qed
\end{enumerate}

\subsection{The specific construction}

We are now ready for our definition of $d$.  Given $u,v \in E$,
suppose that $u = va$ or $u = vb$. We define
\[d(u,v) = 2^{-p_0(u)}.\]

There is one thing we need to check at once.

\begin{lem}
With the above definition of $d$, $S_0$ and $T_0$ are Lipschitz-1
for the metric $\rho$.
\end{lem}

\noindent \textbf{Proof} \hspace{3pt}  We give the argument for
$S_0$; the argument for $T_0$ is similar. Let $u,v \in \G_0$. We
see at once from the definition of $S_0$ that $S_0(u) \wedge
S_0(v) = S_0(u \wedge v)$, and so we have
\[\rho(S_0(u),S_0(v)) = \rho(S_0(u),S_0(u \wedge v)) + \rho(S_0(u \wedge
v),S_0(v)).\] It therefore suffices to treat the case $u\ |\ v$.
Suppose $v = uc_1c_2\cdots c_k$, and write $u_i = v|_{l(u)+i}$ for
$i \leq k$, so $u_0 = u$, $u_k = v$. Then we have
\[\rho(u,v) = \sum_{i=0}^kd(u_i,u_{i+1}) = \sum_{i=0}^k 2^{-p_0(u_i)},\]
and $S_0(u_i) = au_i$ for each $i$ so similarly
\begin{eqnarray*}
\rho(S_0(u),S_0(v)) &=& \sum_{i=0}^kd(S_0(u_i),S_0(u_{i+1})) =
\sum_{i=0}^k 2^{-p_0(S_0(u_i))}\\ &=& \sum_{i=0}^k 2^{-p_0(au_i)}
\leq \sum_{i=0}^k 2^{-p_0(u_i)} = \rho(u,v),
\end{eqnarray*}
using Lemma \ref{lem:perprop}.  This gives the required result.
\qed

\vspace{10pt}

\noindent \textbf{Remark} \hspace{3pt} It follows that $S_0$,
$T_0$ can be extended to the completion $\G$.

In addition, $\L_d$ contains no periodic infinite words, for if
$z$ is a periodic infinite word with repeat block $r$ then
\[l_d(z) = \sum_{i \geq 0}d(z|_i,z|_{i+1}),\]
with all terms equal to $2^{-l(r)}$ for $i \geq l(r)$, so that
$l_d(z) = \infty$.

Hence we have Lipschitz-1 maps $S$,$T$ on the complete metric
space $(\G,\rho)$, no word in which has a fixed point (since $\G$
contains no periodic infinite words).  Finally we come to the
punchline.

\begin{thm}
For any $x,y \in \G$ we have
\[\min\{\rho(S(x),S(y)), \rho(T(x),T(y))\} \leq \frac{3}{4}\rho(x,y).\]
\end{thm}

\noindent \textbf{Proof} \hspace{3pt}  If $x=y$ this is trivial,
so suppose $x \not= y$.  Then $x \wedge y$ is a finite word,
irrespective of the lengths of $x$ and $y$, and we have
\begin{eqnarray*}
\rho(x,y) &=& \rho(x,x\wedge y) + \rho(x \wedge y,y)\\ &=&
\sum_{i=l(x \wedge y)}^{l(x)-1}2^{-p_0(x|_i)} + \sum_{i=l(x \wedge
y)}^{l(y)-1}2^{-p_0(y|_i)}\\ &=& \sum_{e \in A}d(e),
\end{eqnarray*}
where $A = A_1 \cup A_2$ is the union of the set $A_1$ of
individual edges on the path from $x\wedge y$ to $x$ and the set
$A_2$ of edges on the path from $x \wedge y$ to $y$. By Part 2 of
Lemma \ref{lem:perprop} we may partition $A$ into two subsets
$B_a$ and $B_b$, where $e = \{u,v\}$ with $u \leq v$ is in $B_a$
(respectively, in $B_b$) if $p_0(u) = p_0(au)$ (respectively, if
$p_0(u) = p_0(bu)$).  Crucially this implies that if $\{u,v\} \in
B_b$ then $p_0(au) \geq p_0(u) + 1$ and so
\[d(S(e)) = 2^{-p_0(S(u))} \leq \frac{1}{2}2^{-p_0(u)} =
\frac{1}{2}d(e)\] (writing $S(e)$ for the edge $\{S(u),S(v)\}$),
and conversely for $\{u,v\} \in B_a$. Since $S,T$ are Lipschitz-1
we always have $d(S(e)),d(T(e)) \leq d(e)$ for any $e$.

Now we have two cases.  If
\[\sum_{e \in B_a}d(e) \geq \sum_{e \in B_b}d(e)\]
then
\begin{eqnarray*}
\rho(T(x),T(y)) &=& \sum_{e \in B_a} d(T(e)) + \sum_{e \in B_b}
d(T(e)) \leq \frac{1}{2}\sum_{e \in B_a}d(e) + \sum_{e \in
B_b}d(e) \leq \frac{3}{4}\rho(x,y),
\end{eqnarray*}
as required. On the other hand if
\[\sum_{e \in B_a}d(e) \leq \sum_{e \in B_b}d(e)\]
then the above inequality holds with $S$ in place of $T$.  This
completes the proof. \qed

\section{The case of commuting self-maps of a complete metric
space}\label{sec:commcaseprog}

We now consider the problem with finitely many commuting self-maps
$T_1,\ldots,T_n$ of a complete metric space $X$.  In \cite{Stein}
it is shown that some composition of the maps $T_1,...,T_n$ has a
fixed point if these maps are all Lipschitz and their Lipschitz
constants satisfy certain conditions.

The approach we take here is to ask for more: we set out to find a
common fixed point for all the maps in $\F$.  This is motivated by
the following very simple result for the case of $X$ compact:

\begin{prop}
Let $(X,\rho)$ be a compact metric space and $\F$ a (not
necessarily finite) commuting family of continuous self-maps of
$X$ such that for any $x,y \in X$ there is some $U \in F$ with
$\rho(U(x),U(y)) \leq \rho(x,y)$.  Then exactly one point of $X$
is fixed by every member of $U$.
\end{prop}

\noindent \textbf{Proof} \hspace{3pt} The proof of this is fairly
routine.

First we show by contradiction that any member of $\F$ has a fixed
point.  To this end suppose $T \in \F$ has no fixed point.  Since
$T$ is continuous, so is the function $f:X \to [0,\infty): x
\mapsto \rho(x,T(x))$, and since $X$ is compact $f$ must attain
its minimum, say at $x \in X$ with $\rho(x,T(x)) = \alpha$. Since
$T$ has no fixed point, $\alpha > 0$, but by hypothesis we may
choose $S \in \F$ with
\[\rho(S(x),T(S(x))) = \rho(S(x),S(T(x))) < \rho(x,T(x)) = \alpha,\]
contradicting the definition of $\alpha$.  So $T$ must have a
fixed point.

Next we observe that if $T \in \F$ and $x$ is a fixed point for
$T$, then so is $S(x)$ for any $S \in \F$, since members of $\F$
commute.  Given a finite subset $\cal{G}$ of $\F$, let
$A_{\cal{G}}$ be the set of common fixed points in $X$ of all $T
\in \cal{G}$. By the arguments above we know that $A_{\{T\}} \not=
\emptyset$ for every $T \in \F$ and that $S(A_{\{T\}}) \subseteq
A_{\{T\}}$ for every $S,T \in \F$. Furthermore $A_{\{T\}}$ is
closed, hence also compact, and so the first part of the proof
applied with $A_{\{T\}}$ in place of $X$ shows that
$S|_{A_{\{T\}}}$, and hence $S$, must have a fixed point in
$A_{\{T\}}$. Therefore $A_{{\{S,T\}}} \not= \emptyset$. An
induction on $|\cal{G}|$ now gives that $A_{\cal{G}} \not=
\emptyset$ for any finite subset $\cal{G}$ of $\F$.

Finally we observe that $(A_{\cal{G}})_{\cal{G} \in [\F]^{<
\omega}}$ is a family of compact sets with the finite intersection
property, hence has non-empty overall intersection, say $A$.  Any
member of $A$ is then a fixed point for every member of $\F$, and
as in the original Banach contraction principle we see that $A$
must be a singleton. \qed

\vspace{10pt}

The above proof relies on the compactness hypothesis in several
places, and so we cannot hope to gain much by modifying it.  In
the sequel we construct a different argument to show that in the
general commuting case the result holds for $n=2$.  It is worth
noting that once we drop compactness the restriction to a finite
family really is needed, as shown by the example in \cite{Stein},
p.736, of a self-map $T$ of $[1,\infty)$ such that for any $\g >
0$ and $x,y \in [1,\infty)$ there is some $n \geq 1$ with $|T^n(x)
- T^n(y)| \leq \g|x-y|$, and yet such that no $T^k$ has a fixed
point.

We start by giving a convenient necessary and sufficient condition
for the existence of a common fixed point, and then show that it
is satisfied for continuous self-maps if $\F$ has size 2. We aslo
show that that the result holds if all the members of $\F$ are
uniformly continuous and the space $X$ is bounded.

\subsection{A necessary and sufficient condition for a common fixed point}

The following simple lemma turns out to be a crucial step.

\begin{lem}\label{lem:smalltriang}
Suppose $\F = \{T_1,\ldots,T_n\}$ is a finite $\g$-contractive
family on a complete metric space $(X,\rho)$ such that for any
$\varepsilon > 0$ there is some $x \in X$ with $\rho(x,T_i(x))
\leq \varepsilon$ for each $i \leq n$.  Then $\F$ has a common
fixed point.
\end{lem}

\noindent \textbf{Remark} \hspace{3pt} The last conditions in the
lemma above can be described by saying that $\F$ has arbitrarily
good ``approximate common fixed points''.  We will call $x \in X$
an \textbf{$\varepsilon$-approximate common fixed point} for $\F$
if each member of $\F$ moves $x$ by at most $\varepsilon$, as
above. This lemma tells us that the existence of arbitrarily good
approximate common fixed points implies the existence of a common
fixed point. The proof is pleasantly straightforward.

\vspace{10pt}

\noindent \textbf{Proof} \hspace{3pt} For each $m \geq 0$ choose
$x_m \in X$ such that $\rho(x_m,T_i(x_m)) \leq 2^{-m}$ for each $i
\leq n$.  We will show that $(x_m)_{m \geq 0}$ is a Cauchy
sequence (without any further assumptions on our choice of $x_m$).
Suppose $p > q \geq 0$, and consider $\rho(x_p,x_q)$. Since $\F$
is $\g$-contractive we know that we may choose $i \leq n$ with
$\rho(T_i(x_p),T_i(x_q)) \leq \g \rho(x_p,x_q)$.  Combining this
with the triangle inequality gives
\begin{eqnarray*}
\rho(x_p,x_q) &\leq& \rho(x_p,T_i(x_p)) + \rho(T_i(x_p),T_i(x_q))
+ \rho(T_i(x_q),x_q)\\ &\leq& \g \rho(x_p,x_q) + 2^{-p} + 2^{-q},
\end{eqnarray*}
hence
\[(1-\g)\rho(x_p,x_q) \leq 2^{-p} + 2^{-q} < 2^{-(q-1)},\]
and so $\rho(x_p,x_q) \leq \frac{1}{1-\g}2^{-(q-1)}$.  This may be
made arbitrarily small by choosing $q$ sufficiently large, and so
our sequence is Cauchy, as required.  Since $T_1,\ldots,T_n$ are
continuous it follows at once from the definition of an
approximate common fixed point that the limit is the desired
common fixed point, so we are done. \qed

\subsection{Two continuous commuting maps}

Let our two self-maps be $S$ and $T$.  We start with a few
definitions, the significance of which will become clear later.

\begin{dfn}
A point $x \in X$ is said to be a \textbf{sidestepping point} for
$T$ if $\rho(T(x),T^2(x)) > \g \rho(x,T(x))$.  A sidestepping
point of $S$ is defined analogously.
\end{dfn}

That is, $x$ is a sidestepping point for $T$ if $T$ does
\emph{not} contract the pair $(x,T(x))$ by a factor of $\g$; in
this case the contractive assumption shows that
$\rho(S(x),S(T(s))) = \rho(S(x),T(S(x))) \leq \g \rho(x,T(x))$.

We will prove the following result by showing that either we may
find a fixed point for $T$ directly or the conditions of Lemma
\ref{lem:smalltriang} are satisfied. This will require a careful
estimate.

\begin{thm}\label{thm:mainfortwo}
Suppose $\F = \{S,T\}$ is a commuting contractive family on
$(X,\rho)$.  Then precisely one point of $X$ is fixed by both $S$
and $T$.
\end{thm}

This result follows from the following.

\begin{prop}\label{prop:mainfortwo}
Suppose $\F = \{S,T\}$ is a commuting contractive family on
$(X,\rho)$.  Then:
\begin{enumerate}
\item either there is some $x \in X$ such that $(T^i(x))_{i \geq
0}$ is a Cauchy sequence, or;
\item $X$ contains an $\varepsilon$-approximate common fixed point
for $\F$ for every $\varepsilon > 0$.
\end{enumerate}
(It is possible that both of these hold).
\end{prop}

\noindent \textbf{Proof} \hspace{3pt} Suppose that $(T^i(x))_{i
\geq 0}$ is not Cauchy for any $x \in X$. In particular, this
implies that for any $x \in X$, some $T^i(x)$ is a sidestepping
point for $T$; for otherwise we would have
$\rho(T^{i+1}(x),T^{i+2}(x)) \leq \g\rho(T^i(x),T^{i+1}(x))$ for
every $i \geq 0$, easily implying that $(T^i(x))_{i \geq 0}$ is
Cauchy.

We define the sequence $(x_n)_{n \geq 0}$ recursively as follows.
Set $x_0 = x$.  Suppose we have defined $x_0,x_1,\ldots,x_{n-1}$.
If $x_{n-1}$ is not a sidestepping point for $T$, let $x_n =
T(x_{n-1})$; otherwise let $x_n = S(x_{n-1})$. See the sketch
below. This definition guarantees that $\rho(x_n,T(x_n)) \leq \g
\rho(x_{n-1},T(x_{n-1}))$ for every $n \geq 1$. It follows that
$\rho(x_n,T(x_n)) \to 0$ as $n \to \infty$, and so for any $\eta >
0$ we may assume that $\rho(x,T(x)) < \eta$ by starting at $x_n$
for some sufficiently large $n \geq 0$, rather than at $x$.

Now, let $n_1 \leq n_2 \leq n_3 \leq \ldots$ be those $n$ which
index sidestepping points for $T$ in the sequence $(x_n)_{n \geq
0}$.  We show that if $\rho(x,T(x)) \leq \eta$, then there is a
sidestepping points $x_{n_k}$ of our sequence with
\[\rho(x_{n_k},S(x_{n_k})) \leq \frac{3}{(1 - \g)^2}\eta +
\frac{1}{1-\g}\eta.\] Thus, from our sequence $(x_n)_{n \geq 0}$,
constructed to put a tight bound on $\rho(x_n,T(x_n))$, we may
choose points $x_n$ for which $\rho(x_n,S(x_n))$ is also small;
these are our approximate common fixed points.

We will obtain this bound by considering the distances between
successive sidestepping points in the sequence $(x_n)_{n \geq 0}$;
that is, the distances $\rho(x_{n_k},x_{n_{k+1}})$ for $k \geq 1$.
We first prove that there exist $k \geq 1$ for which this distance
is at most $\frac{3}{(1-\g)^2}\eta$. Set $L_k =
\rho(x_{n_k},x_{n_{k+1}})$ for $k \geq 1$. Since $\F$ is a
contractive family, one of $S$, $T$ contracts the distance
$\rho(x_{n_k},x_{n_{k+1}})$.  There are two cases.

\begin{enumerate}
\item Suppose first that for some $k \geq 1$, $T$ contracts the distance
$\rho(x_{n_k},x_{n_{k+1}})$. In this case we have the picture
shown in Fig 2.

\begin{figure}
\begin{center}
\includegraphics[width = 0.8\textwidth]{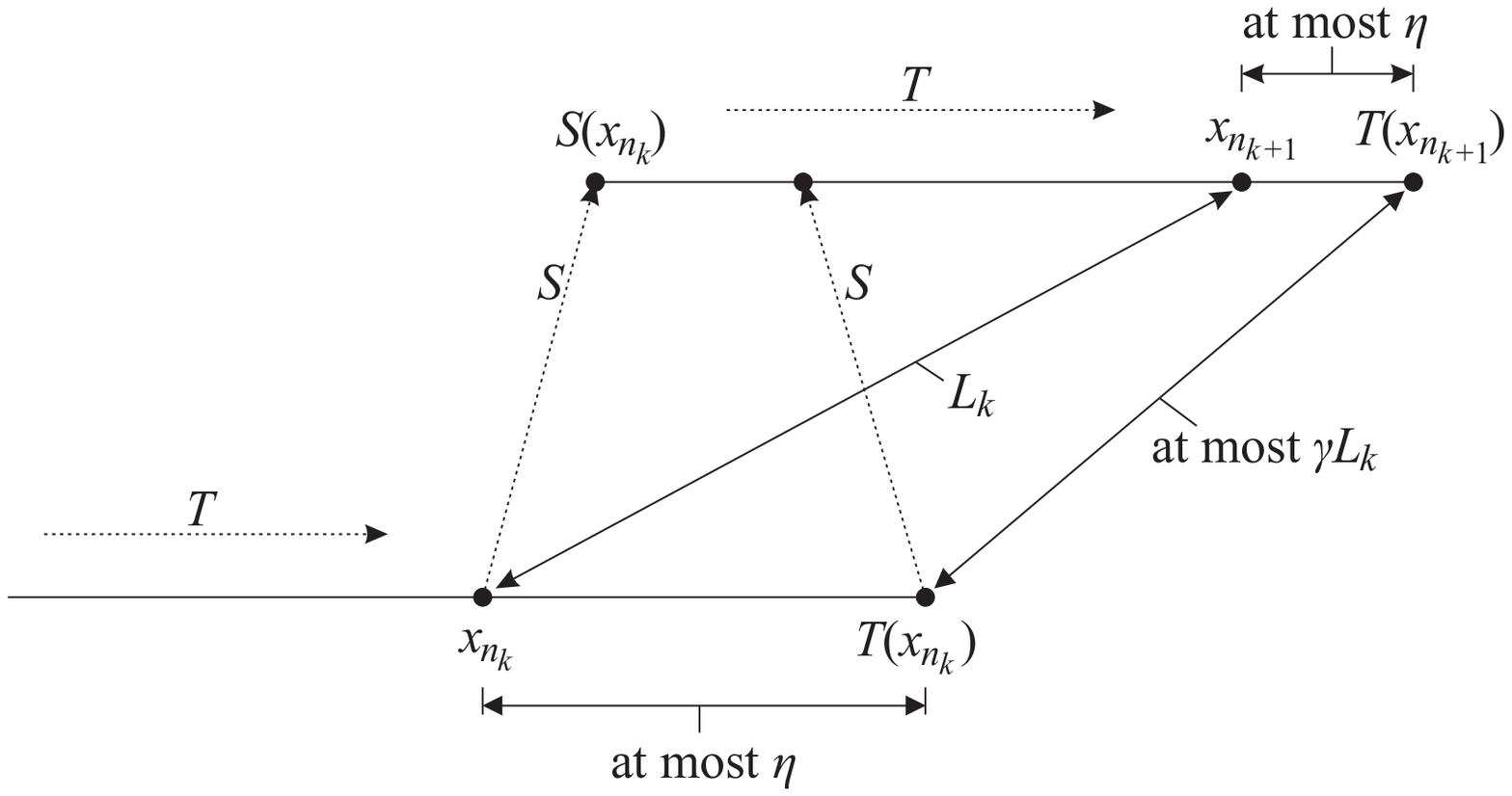}
\caption{Case 1}
\end{center}
\end{figure}

The triangle inequality gives
\begin{eqnarray*}
L_k = \rho(x_{n_k},x_{n_{k+1}}) &\leq& \rho(x_{n_k},T(x_{n_k})) +
\rho(T(x_{n_k}),T(x_{n_{k+1}})) +
\rho(T(x_{n_{k+1}}),x_{n_{k+1}})\\ &\leq& \eta + \g L_k + \eta,
\end{eqnarray*}
and so $L_k$ must itself be at most $\frac{2}{1-\g}\eta \leq
\frac{3}{(1 - \g)^2}\eta$, and we are home;

\item Alternatively, suppose that $T$ does not contract
$\rho(x_{n_k},x_{n_{k+1}})$ for any $k$, so they must all be
contracted by $S$. See Fig 3.

\begin{figure}
\begin{center}
\includegraphics[width = 0.8\textwidth]{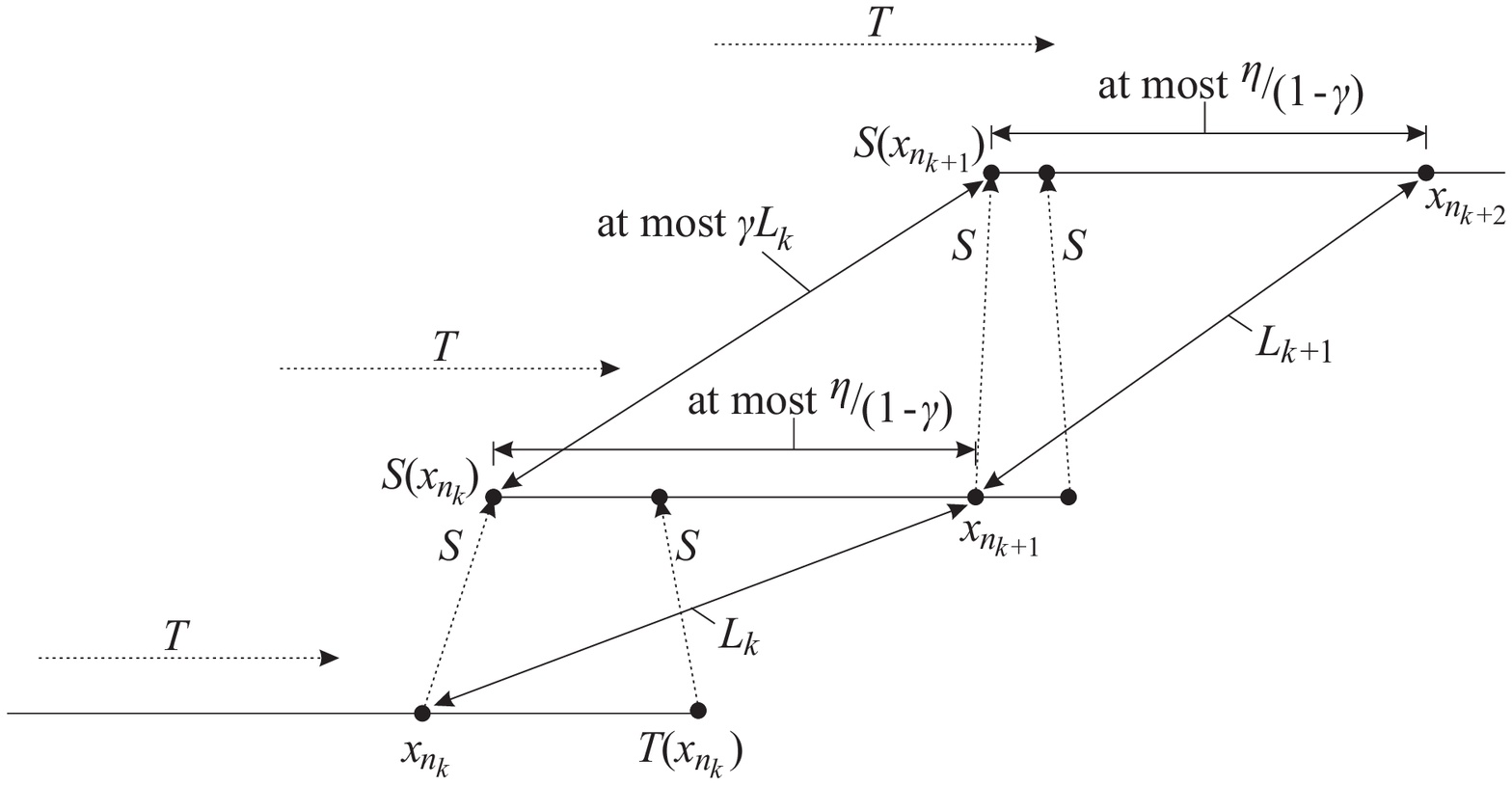}
\caption{Case 2}
\end{center}
\end{figure}

This time the triangle inequality gives the following estimate:
\begin{eqnarray*}
\rho(x_{n_{k+1}},x_{n_{k+2}}) &\leq& \rho(x_{n_{k+1}},S(x_{n_k}))
+ \rho(S(x_{n_k}),S(x_{n_{k+1}})) +
\rho(S(x_{n_{k+1}}),x_{n_{k+2}})\\ &=& \rho(x_{n_{k+1}},x_{n_k+1})
+ \rho(S(x_{n_k}),S(x_{n_{k+1}})) +
\rho(x_{n_{k+1}+1},x_{n_{k+2}})\\ &\leq& \sum_{n_k+1 \leq i <
n_{k+1}}\rho(x_i,x_{i+1}) + \g L_k + \sum_{n_{k+1}+1 \leq i <
n_{k+2}}\rho(x_i,x_{i+1})\\ &\leq& \g L_k + \frac{2}{1-\g}\eta.
\end{eqnarray*}
Hence $L_{k+1} = \rho(x_{n_{k+1}},x_{n_{k+2}}) \leq \g L_k +
\frac{2}{1-\g}\eta$. Similarly, $L_{k + 2} \leq \g L_{k+1} +
\frac{2}{1-\g}\eta$, and so on.  It is now easy to see that if
$L_k > \frac{3}{(1-\g)^2}\eta$ then we obtain a finite decreasing
sequence $L_k,L_{k+1},\ldots$ until we reach an $i$ with $L_{k+i}
\leq \frac{3}{(1-\g)^2}\eta$, and so are done.
\end{enumerate}
Having found a suitable $x_{n_k}$ by the above procedure, it
remains only to observe that selecting $\eta$ so small that both
$\eta < \varepsilon$ and
\[\Big(\frac{1}{1-\g} + \frac{3}{(1-\g)^2}\Big)\eta < \varepsilon,\]
we have
\[\rho(x_{n_i},S(x_{n_i})),\rho(x_{n_i},T(x_{n_i})) <
\varepsilon.\] Thus $x_{n_i}$ is an $\varepsilon$-approximate
common fixed point for $\F$. This completes the proof of the
proposition. \qed

\vspace{10pt}

\noindent \textbf{Proof of Theorem
\ref{thm:mainfortwo}}\hspace{3pt} We have two cases, corresponding
to the dichotomy of Proposition \ref{prop:mainfortwo}.

Suppose first that there is some $x \in X$ for which $(T^i(x))_{i
\geq 0}$ is Cauchy, with limit $y$. Then the continuity of $T$
gives at once that $T(y) = y$. Let $A$ be the set of all fixed
points of $T$ in $X$, so $A$ is nonempty. Since $S,T$ commute,
$S(A) \subseteq A$. In addition, $T|_A = \rm{Id}_A$, so no two
points of $A$ are contracted by $T$, and therefore any such pair
of points must be contracted by $S$. Hence $S|_A$ is a contraction
of $A$, and so has a unique fixed point in $A$, by the original
contraction mapping principle.  This is then the unique common
fixed point for $S$ and $T$.

If we cannot choose such an $x \in X$, then Proposition
\ref{prop:mainfortwo} yields $\varepsilon$-approximate common
fixed points for $\F$ for every $\varepsilon > 0$, and so the
desired existence follows from Lemma \ref{lem:smalltriang}.
Uniqueness is proved as in Banach's original result. \qed

\vspace{10pt}

This completes the proof of our main positive result, that an
contractive family $\F$ of size two has precisely one common fixed
point.  This leaves us with the following.

\begin{conj}
If $\F$ is any finite commuting contractive family of self-maps on
a complete metric space $X$ then $\F$ has exactly one common fixed
point.
\end{conj}

\subsection{Other conditions}

In this last subsection we will prove the existence of a common
fixed point for a commuting contractive family $\F$ of arbitrary
finite size under the two additional assumptions that the members
of $\F$ are uniformly continuous and that $X$ is bounded. Once
again the proof rests on Lemma \ref{lem:smalltriang}.

\begin{lem}
Suppose $X$ is bounded and $\F$ is a finite commuting contractive
family of uniformly continuous self-maps of $X$. Then $X$ contains
a $\varepsilon$-approximate common fixed point for $\F$ for any
$\varepsilon > 0$.
\end{lem}

\noindent \textbf{Proof} \hspace{3pt} Let $\F =
\{T_1,T_2\,\ldots,T_n\}$. We shall prove by induction on $k$ that
for every $k \leq n$ there is some $x \in X$ with
\[\rho(x,T_1(x)),\rho(x,T_2(x)),\ldots,\rho(x,T_k(x)) < \varepsilon;\]
the case $k = n$ is then the desired conclusion.

We have seen in previous arguments (such as in Proposition
\ref{prop:mainfortwo}) that there are point $x$ for which
$\rho(x,T_1(x))$ is arbitrarily small; this is the base case
$k=1$. Now assume that the result holds for some $k \leq n-1$. Let
$M < \infty$ be the diameter of $X$, and choose $m \geq 0$ so that
$\g^mM < \varepsilon$. Since the members of $\F$ are uniformly
continuous, there is some $\delta > 0$ such that
\[\rho(U(x),U(y)) < \varepsilon\]
\emph{whenever} $\rho(x,y) < \delta$ and $U$ is word of length at
most $m$ in the members of $\F$. By the inductive hypothesis, we
may choose a $\delta$-approximate common fixed point for
$\{T_1,T_2,\ldots,T_k\}$. By the definition of $M$, we have
$\rho(x,T_{k+1}(x)) \leq M$. By applying an appropriately-chosen
member of $\F$, say $S$, we may obtain $\rho(S(x),S(T_{k+1}(x)))
\leq \g M$. Continuing in this way, after at most $m$ steps we
will have a word $U$ of length at most $m$ in the members of $\F$
such that $\rho(U(x),U(T_{k+1}(x))) \leq \g^mM < \varepsilon$. By
our choice of $\delta$, we know that we also have
$\rho(U(x),U(T_{j}(x))) < \varepsilon$ for $j = 1,2,\ldots,k$.
Hence the point $U(x)$ witnesses the desired result for $k+1$, and
the induction continues. \qed

\vspace{10pt}

Arguing uniqueness in the same way as always, we deduce from the
above:

\begin{cor}
Suppose $X$ is bounded and $\F$ is a finite commuting contractive
family of uniformly continuous self-maps of $X$. Then $\F$ has a
unique common fixed point in $X$. \qed
\end{cor}

\section*{Acknowledgements}

The above work was carried out under a summer research studentship
funded by Trinity College, Cambridge over the long vacation period
of 2004.  My thanks go to Dr I. Leader (Department of Pure
Mathematics and Mathematical Statistics, University of Cambridge)
for the supervision of the project and the many protracted
discussions this entailed, and to Trinity College for their
support over this long vacation period.

My thanks also go to David Fremlin for his clear explanation of
the Generalized Banach Contraction Theorem.


\begin{flushright}
\emph{Tim D. Austin}

\emph{Trinity College}

\emph{Cambridge}

\emph{August 2004}
\end{flushright}

\end{document}